\documentclass{amsart}
\usepackage{amsmath}
\usepackage{amsthm}
\usepackage{amsfonts}
\usepackage{amssymb}
\usepackage{latexsym}
\usepackage{enumerate}

\newtheorem{theo}{Theorem}[section]

\newtheorem{lem}[theo]{Lemma}

\begin{document}

\title[The Zeckendorf sum of digits of powers]{Extremal orders of the Zeckendorf sum of digits of powers}
\author{Thomas Stoll}
\address{Institut de Math\'ematiques de Luminy, Universit\'e de la M\'editerran\'ee, 13288 Marseille Cedex 9, France,}
\email{stoll@iml.univ-mrs.fr}


\maketitle

\begin{abstract}
  Denote by $s_F(n)$ the minimal number of Fibonacci numbers needed to write $n$ as a sum of Fibonacci numbers.
  We obtain the extremal minimal and maximal orders of magnitude of $s_F(n^h)/s_F(n)$ for any $h\geq 2$. We use this
  to show that for all $N>N_0(h)$ there is a $n$ such that $n$ is the sum of $N$ Fibonacci numbers and $n^h$ is the sum
  of at most $130 h^2$ Fibonacci numbers. Moreover, we give upper and lower bounds on the number of $n$'s with small and
  large values $s_F(n^h)/s_F(n)$. This extends a problem of Stolarsky to the Zeckendorf representation of powers, and 
  it is in line with the classical investigation of finding perfect powers among the Fibonacci numbers and their finite sums.
\end{abstract}

\section{Introduction}

Denote by $s_q(n)$ the sum of digits in the usual $q$-ary digital expansion of $n$. Stolarsky~\cite{St78} studied the maximal 
and minimal order of magnitude of the ratio $s_2(n^h)/s_2(n)$ for fixed $h\geq 2$. It is reasonable to expect that the quantities $s_2(n^h)$ and 
$s_2(n)$ are independent in the sense that the $\limsup$ of the ratio tends to $\infty$ and the $\liminf$ to $0$ as $n$ tends 
to infinity. It is an interesting question to find the extremal orders of magnitude of this ratio. In a recent work, Hare, Laishram and the author~\cite{HLS10} were able to settle an open question of Stolarsky, so to finally get a complete picture of the \textit{maximal and minimal order of magnitude} of the ratio $s_q(n^h)/s_q(n)$.

\begin{theo}[\cite{St78, HLS10}]\label{theoA}
  There exist $c_1$ and $c_2$, depending at most on $q$ and $h$, such that for all $n\geq 2$,
    $$\frac{c_2}{\log n}\leq \frac{s_q(n^h)}{s_q(n)}\leq c_1 (\log n)^{1-1/h}.$$
  This is best possible in that there exist $c_1'$ and $c_2'$, depending at most on $q$ and $h$, such that
  \begin{equation}\label{c1dot}
    \frac{s_q(n^h)}{s_q(n)}> c_1' (\log n)^{1-1/h},
  \end{equation}
  respectively,
  \begin{equation}\label{c2dot}
    \frac{s_q(n^h)}{s_q(n)}< \frac{c_2'}{\log n}
  \end{equation}
  infinitely often.
\end{theo}

In the present paper we find the maximal and minimal order of magnitude of the ratio $s_F(n^h)/s_F(n)$, where $s_F$ denotes
the Zeckendorf sum of digits function, and we give a Diophantine application.
Let
\begin{equation}\label{xdef}
  x=\sum_{2\leq j\leq n} \varepsilon_j F_j,
\end{equation}
with $\varepsilon_n=1$ and $\varepsilon_j\in \{0,1\}$ be the (greedy) Zeckendorf expansion of $x\in \mathbb{Z}^+$ with respect to the Fibonacci numbers $F_j$.
Recall that in this expansion we do not allow adjacent $1$ digits~\cite{Ze72, AS03}. Hence $x$ can have at most $\lfloor n/2 \rfloor$ digits $1$ in its
expansion. We write $x=(\varepsilon_n \varepsilon_{n-1}\ldots \varepsilon_2)_F$ to refer to this expansion. Denote by $s_F$ the Zeckendorf sum of
digits function defined by $$s_F(x)=\sum_{2\leq j\leq n}\varepsilon_j.$$
This function can also be interpreted as the minimal number of Fibonacci numbers needed to write $n$ as a sum of Fibonacci numbers.
$s_F$ shares many properties with the ordinary sum of digits function $s_q$. For instance, $s_F$ is also subadditive (i.e., $s_F(a+b)\leq s_F(a)+s_F(b)$ for all $a,b\geq 1$),
has fractal summatory behaviour~\cite{CV86} and satisfies a Newman phenomenon~\cite{DS00}.  Contrary to $s_q$~\cite{HLS10}, the function $s_F$ is not submultiplicative, as the example
$$ 2\cdot 3 =(10)_F \cdot (100)_F = (1001)_F=6$$ shows. 
Therefore, there is \textit{a priori} no obvious relation between $s_F(n^h)$ and $s_F(n)^h$. Drmota and Steiner~\cite{DS02}, extending a result
of Bassily and K\'atai~\cite{BK95}, showed that $s_F(n^h)$ properly renormalized
is asymptotically normally distributed. The mean value of $s_F(n^h)$ is asymptotically $h$ times the mean value of $s_F(n)$ which is $c_F \log n$ with a suitable constant $c_F$~\cite{PT89}.
This means that we expect $n^h$ to have roughly $h$ times as many 1's in the Zeckendorf expansion compared to $n$, thus
the ratio $s_F(n^h)/s_F(n)$ should be roughly $h$. 
Our main result is as follows.

\begin{theo}\label{theom1}
  There exist $c_3$ and $c_4$, depending at most on $h$, such that for all $n\geq 2$,
\begin{equation}\label{c3}
  \frac{c_4}{\log n} \leq \frac{s_F(n^h)}{s_F(n)}\leq c_3 \log n.
\end{equation}
  This is best possible in that there exist $c_3'$ and $c_4'$, depending at most on $h$, such that
\begin{equation}\label{c3dot}
  \frac{s_F(n^h)}{s_F(n)}> c_3' \log n
\end{equation}
respectively,
\begin{equation}\label{c4dot}
  \frac{s_F(n^h)}{s_F(n)}<\frac{c_4'}{\log n},
\end{equation}
  infinitely often. Moreover, possible values for the constants are
\begin{equation}\label{choose}
  c_3=2h, \qquad c_3'=1, \qquad c_4=\frac{1}{2}, \qquad c_4'=120 h^2.
\end{equation}
\end{theo}

This is strongly related to the classical investigation of finding perfect powers among Fibonacci numbers and their finite
sums. A recent deep result of Bugeaud, Mignotte and Siksek~\cite{BMS06} says that the only powers
$n^h$ that are Fibonacci numbers (or equivalently, with $s_F(n^h)=1$), are $1,8$ and $144$. 
From~(\ref{c4dot}),~(\ref{choose}) and our construction we obtain the following Diophantine result.
\begin{theo}\label{fibcoro}
 For any $h\geq 2$ there exists $N_0(h)$, only depending on $h$, such that for all $N>N_0$ there exists an integer $n$ with the following two properties:
 \begin{enumerate}
   \item[(i)] $n$ is the sum of $N$ distinct, non-adjacent Fibonacci numbers.
   \item[(ii)] $n^h$ is the sum of \emph{at most} $130 h^2$ Fibonacci numbers.
 \end{enumerate}
\end{theo}

Recently, Bugeaud, Luca, Mignotte and Siksek~\cite{BLMS08} found
all powers which are at most one away from a Fibonacci number. 
In our context, this is the investigation of finding powers
with \textit{very large} and \textit{very small} sum of digits values. 
A refinement of our construction yields that $s_F(n^h)$ is
\textit{small} and \textit{large} indeed quite often compared to $s_F(n)$. 

\begin{theo}\label{theom2}
  For $\varepsilon>0$ there exists $$\alpha > \frac{1}{\max\left(36h^2/\varepsilon+18, 8h+1\right)}$$ such that
  \begin{equation}\label{small}
    \#\{n<N:\quad \frac{s_F(n^h)}{s_F(n)} < \varepsilon\}\gg N^\alpha.
  \end{equation}
\end{theo}

\begin{theo}\label{theom3}
  For $\delta>0$ there exists $$\beta > \frac{1}{h(\delta+1)+2}$$ such that
  \begin{equation}\label{large}
    \#\{n<N:\quad \frac{s_F(n^h)}{s_F(n)} > \delta\}\gg N^{\beta}.
  \end{equation}
\end{theo}

In Section~\ref{secprelim} we collect and state some facts about Fibonacci numbers, Lucas numbers and the Zeckendorf sum
of digits function, which we will need in the proofs. In Sections~\ref{sectheom1} and~\ref{sectheom2} we then give the elementary constructions that
prove~(\ref{c3dot}),~(\ref{c4dot}) and Theorem~\ref{fibcoro}. Section~\ref{sectheom3} is devoted to the proofs of Theorems~\ref{theom2} and~\ref{theom3}.

\section{Preliminaries}\label{secprelim}

Since $F_n=\lfloor \phi^n/\sqrt{5}\rfloor$, where $\phi =\frac{1}{2}(\sqrt{5}+1)$ is the golden ratio,
we have by~(\ref{xdef}) that
$$\frac{\phi^{n-3/2}}{\sqrt{5}}< \left \lfloor \frac{\phi^n}{\sqrt{5}}\right \rfloor \leq x <
\left \lfloor \frac{\phi^{n+1}}{\sqrt{5}}\right \rfloor\leq \frac{\phi^{n+1}}{\sqrt{5}}$$
for $n\geq 2$. Therefore,
\begin{equation}\label{numdig}
  n=\frac{\log x}{\log \phi}+\gamma_n,
\end{equation}
where $\gamma_n$ lies in the interval
$$(\delta, \delta'):=\left(\frac{\log \sqrt{5}}{\log \phi}-1, \frac{\log \sqrt{5}}{\log \phi}+\frac{3}{2}\right)\approx (0.672, 3.172).$$ This already implies~(\ref{c3})
with $c_3=2h$ and $c_4= \frac{1}{2}$. 

In the following we show that substracting a ``small'' number from a Fibonacci number
gives rise to a large number of digits $1$ in the Zeckendorf expansion.

\begin{lem}\label{expand}
  Let $k\geq 1$.
  \begin{enumerate} \item[(i)] For $0<z\leq F_{2k+1}$ we have
  $$s_F(F_{2k+1}-z)=k-l+s_F(F_{2l+1}-z)\geq k-\frac{\log z}{2 \log \phi}-\frac{\delta'}{2},$$
  where $l$ is such that $F_{2l}<z\leq F_{2l+1}$.
  \item[(ii)] For $0<z\leq F_{2k}$ we have
  $$s_F(F_{2k}-z)=k-l+s_F(F_{2l}-z)\geq k-\frac{\log z}{2 \log \phi}-\frac{\delta'}{2},$$
  where $l$ is such that $F_{2l-1}<z\leq F_{2l}$.
  \end{enumerate}
\end{lem}
\begin{proof}
  Part $(i)$ follows at once from the identity
  $$F_{2k+1}-z=\left(\sum_{i=l+1}^k F_{2i}+F_{2l+1}\right)-z=\sum_{i=l+1}^k F_{2i}+\left(F_{2l+1}-z\right)$$ and~(\ref{numdig}). The second case is similar.
\end{proof}

Denote by $L_k$ the Lucas numbers defined by
\begin{equation}\label{lucasdef}
  L_k=F_{k-1}+F_{k+1}=\lfloor \phi^k \rfloor.
\end{equation}
Powers and products of Lucas numbers are given by the following formul\ae.
\begin{lem}\label{lucasprop}
For all $k> l\geq 1$ and $h\geq 2$ we have
\begin{equation}\label{powerformula}
  L_k^h = \frac{1}{2} \sum_{i=0}^h \binom{h}{i} (-1)^{ik} L_{(h-2i)k},
\end{equation}
\begin{equation}\label{id2}
  L_k L_l=L_{k+l}+(-1)^l L_{k-l}.
\end{equation}
\end{lem}

\begin{proof}
  See for example~\cite{Va08}.
\end{proof}

Formula~(\ref{powerformula}) shows that powers of odd indexed Lucas numbers can be written as linear sum of Lucas numbers having positive coefficients.
Furthermore, from~(\ref{id2}) we have that products of two even indexed Lucas numbers can be rewritten as sums of two single Lucas numbers.
We will further need the fact that fixed multiples of Lucas numbers have bounded sum of digits values.

\begin{lem}\label{lucasmulti}
  For $m\geq 1$ there exists $k_0=k_0(m)$ such that for all $k\geq k_0$,
  $$s_F(mL_k)<\frac{\log m}{\log \phi}+3.$$
\end{lem}

\begin{proof}
  Since $F_l L_k=F_{l+k}-(-1)^l F_{k-l}$ we have that
  \begin{align*}
    F_{2l+1}L_k &=F_{k+2l-1}+F_{k-2l+1},\\
    F_{2l}L_k&=F_{k+2l}-F_{k-2l}=F_{k-2l+1}+F_{k-2l+3}+\cdots+F_{k+2l-1}.
  \end{align*}
  Hence, by writing $m$ in Zeckendorf representation we get that for all $m$ with $F_{2l}< m<F_{2l+1}$ the Zeckendorf representation of $mL_k$
  involves a block of $4l+2$ digits ($k$ sufficiently large) and a following block of zeros only, and for all $m$ with $F_{2l+1}\leq m \leq F_{2l+2}$
  a block of $4l+3$ digits with a block of zeros appended. This yields that for each $m\geq 1$ and $k\geq k_0$ a block of length at most
  \begin{equation}\label{blocklength}
    \frac{2 \log \sqrt{5}m}{\log \phi}+2
  \end{equation}
  appears in the representation. Thus,
  $$s_F(mL_k)\leq \frac{\log \sqrt{5}m}{\log \phi}+1<\frac{\log m}{\log \phi}+3,$$
  which proves the claim.
\end{proof}

\section{Proof of the extremal upper bound}\label{sectheom1}

We use a construction of an extremal sequence based on the power expansion of Lucas numbers~(\ref{powerformula}).
Set $n_k = L_{2k-1}$ for $k\geq 1$. Then by~(\ref{lucasdef}) we have $s_F(n_k)=2$. For the proof of~(\ref{c3dot}) it suffices to show that
$s_F(n_k^h)=2k+O_h(1)$, where the implied constant depends only on $h$.
We have
\begin{align}\label{hexp}
  n_k^h &= F_{h(2k-1)+1}+F_{h(2k-1)-1}\nonumber\\
  &\qquad -\frac{1}{2} \sum_{i=1}^{h-1} \binom{h}{i}(-1)^{(i+1)(2k-1)} L_{(h-2i)(2k-1)}.
\end{align}
The last sum is positive since
\begin{align*}
  \frac{1}{2} \sum_{i=1}^{h-1} \binom{h}{i}(-1)^{(i+1)(2k-1)} &L_{(h-2i)(2k-1)}
  =\left\lfloor\frac{\phi^{h(2k-1)}}{\sqrt{5}}\right\rfloor-\left\lfloor\frac{\phi^{2k-1}}{\sqrt{5}}\right\rfloor^h\\
  &\geq \phi^{h(2k-1)} \left(\frac{1}{\sqrt{5}}-\frac{1}{\sqrt{5}^h}\right)-1>0.
\end{align*}
Moreover, this quantity is small with respect to the leading term. In fact, we get by a trivial estimate
\begin{align*}
  \frac{1}{2} \sum_{i=1}^{h-1} \binom{h}{i}(-1)^{(i+1)(2k-1)} L_{(h-2i)(2k-1)} &\leq 2^{h-1} L_{(h-2)(2k-1)}\\
  &\leq 2^{h-1}\phi^{(h-2)(2k-1)}
\end{align*}
which is smaller than $F_{h(2k-1)-1}$ for sufficiently large $k$. Therefore, using Lemma~\ref{expand}, we get
\begin{align*}
  s_F(n_k^h) &\geq 1+\left\lfloor \frac{h(2k-1)-1}{2} \right\rfloor -\frac{\log (2^{h-1} \phi^{(h-2)(2k-1)})}{2\log \phi}-\frac{\delta'}{2}\\
  &\geq 2k-\frac{h-1}{2}\cdot \frac{\log 2}{\log \phi}-\frac{3}{2}-\frac{\delta'}{2}\\
  &\geq 2k-\frac{3h}{4}-3,
\end{align*}
for $k$ sufficiently large. Therefore, as $k$ tends to infinity,
\begin{align*}
  \frac{s_F(n_k^h)}{s_F(n_k)}&\geq k -\frac{3h}{8}-\frac{3}{2}\\
  &\geq \left(\frac{\log n_k}{2 \log \phi}+\frac{1}{2}\right)-\frac{3h}{8}-\frac{3}{2}\gg \log n_k.
\end{align*}  
Hence, we can put $c_3'=1$ and get~(\ref{c3dot}). \hfill\qed

\section{Proof of the extremal lower bound}\label{sectheom2}

Here, we use a construction which uses~(\ref{id2}).
Let $k\geq 1$ and set
$$n_k= L_{8k}+L_{6k}+L_{4k}+L_{2k}-1.$$
We have
\begin{align}\label{linear}
  s_F(n_k) &=6+s_F(L_{2k}-1)\nonumber\\
  &=6+s_F(F_2+F_4+\cdots+F_{2k-2}+F_{2k+1})=6+k.
\end{align}
First we calculate the Zeckendorf expansion of $n_k^2 = (L_{8k}+L_{6k}+L_{4k}+L_{2k}-1)^2$. We expand
the square by employing~(\ref{id2}) and use the special value $L_0=2$ to get
$$n_k^2 = L_{16k}+2L_{14k}+3L_{12k}+4L_{10k}+L_{8k}+2L_{6k}+3 L_{4k}+4L_{2k}+9.$$
We replace all appearances of multiples of Lucas numbers by the corresponding linear sum in Fibonacci
numbers. In this case, we use
\begin{align*}
  2 L_k&=F_{k+3}+F_{k-3},\\
  3L_k&=F_{k+3}+F_{k+1}+F_{k-1}+F_{k-3},\\
  4L_k&=F_{k+4}+F_{k+1}+F_{k-2}+F_{k-5}.
\end{align*}
It is now a straightforward calculation to write down
the expansion of $n_k^2$. In order to simplify notation, denote by $(e_p e_{p-1} \ldots e_{0})_{l}$ the
sum of Fibonacci numbers $e_p F_{p+l}+e_{p-1} F_{p-1+l}+\cdots+ e_0 F_l$. We get
\begin{align}\label{squares}
  n_k^2 = \quad &(101)_{16k-1} + (1000001)_{14k-3}+(1010101)_{12k-3}\nonumber\\
  &+(1001001001)_{10k-5}+(101)_{8k-1}+(1000001)_{6k-3}\nonumber\\
  &+(1010101)_{4k-3}+(1001001001)_{2k-5}+(10001)_2.
\end{align}
Thus, we have $s_F(n_k^2)=26$ for all $k\geq 7$.

In a similar style we obtain the expansion for $n_k^3$. This time we use~(\ref{id2}) twice to rewrite all products of three Lucas numbers
as sums of four Lucas numbers. We here get
\begin{align*}
  n_k^3 &=L_{24k}+3L_{22k}+6L_{20k}+10L_{18k}+9L_{16k}+9L_{14k}+10L_{12k}\\
  &\quad+12L_{10k}+27L_{8k}+28L_{6k}+27L_{4k}+24L_{2k}+11.
\end{align*}
Similarly as before we replace multiples of Lucas numbers by sums of Fibonacci numbers. We get
\begin{align}\label{cubes}
  n_k^3 = \quad & (101)_{24k-1}+(1010101)_{22k-3}+(10001010001)_{20k-5}\nonumber\\
  &+(10010000001001)_{18k-7}+(10000100101001)_{16k-7}\nonumber\\
  &+(10000100101001)_{14k-7}+(10010000001001)_{12k-7}\nonumber\\
  &+(10100100100001)_{10k-7}+(100100010100001001)_{8k-9}\nonumber\\
  &+(100101000001001001)_{6k-9}+(100100010100001001)_{4k-9}\nonumber\\
  &+(100001010100101001)_{2k-9}+(10100)_2.
\end{align}
For $k\geq 10$ the summands in~(\ref{cubes}) are noninterfering. This yields $s_F(n_k^3)=60$ for $k\geq 10$.
Note that~(\ref{squares}) and~(\ref{cubes}) already prove~(\ref{c4dot}) in the case of $h=2$ and $h=3$.

The general case follows from~(\ref{squares}),~(\ref{cubes}) and Lemma~\ref{lucasmulti}. For that purpose set $h= 2h_1+3h_2$ with $h_1, h_2\geq 0$ 
and consider $n_k^h =(n_k^2)^{h_1}\cdot (n_k^3)^{h_2}$. Since both $n_k^2$ and $n_k^3$ are linear forms in
Lucas numbers with fixed positive coefficients, the powers $(n_k^2)^{h_1}$ and $(n_k^3)^{h_2}$ 
are linear forms with positive coefficients, too, that are independent of $k$. Thus we have that 
$n_k^h$ is a linear form in $4h$ Lucas numbers with positive coefficients independent of $k$ (plus an additive constant).  This means that there exists $k_0=k_0(h)$
 such that for all $k\geq k_0$ the terms in the Lucas sum are noninterfering. All coefficients in this sum are bounded
by $9^h$. Therefore, by Lemma~\ref{lucasmulti},
$$s_F(n_k^h)\leq \left(h\; \frac{\log 9}{\log \phi}+3\right)\cdot(4h+1).$$
Since $\phi^{8k}<n_k \leq \phi^{8k+1}$ we also get
$$s_F(n_k)=k+6 \geq \frac{\log n_k}{8 \log \phi}-\frac{1}{8}+6 \gg \frac{1}{4} \log n_k.$$ 
This shows that for sufficiently large $k$,
$$\frac{s_F(n_k^h)}{s_F(n_k)}<\frac{4 (5h+3)(4h+1)}{\log n_k}<\frac{120 h^2}{\log n_k}.$$
This completes the proof of~(\ref{c4dot}) with $c_4'=120h^2$. \hfill \qed

\medskip

\textit{Proof of Theorem~\ref{fibcoro}:}
This follows at once from~(\ref{linear}) and
$$s_F(n)\leq \frac{\log n}{2\log \phi}+2\ll \frac{13}{12} \log n.$$
\hfill \qed

\section{Proofs of Theorems~\ref{theom2} and~\ref{theom3}}\label{sectheom3}

\textit{Proof of Theorem~\ref{theom2}:} For $m\geq 1$ set $$n_k=n_k(m)=m(L_{8k}+L_{6k}+L_{4k}+L_{2k}-1).$$
As before, we have that $n_k^h$ is a linear sum of Lucas numbers with positive coefficients independent of $k$.
Suppose now
\begin{equation}\label{kcondition}
  k> \frac{h \log(9m)}{\log \phi}+O(1).
\end{equation}
Then the blocks in the expansion of $n_k$ respectively $n_k^h$ are noninterfering. Using~(\ref{blocklength}) we have
$$  s_F(n_k) \geq k-\frac{2 \log m}{\log \phi}+O(1)$$
and
$$ s_F(n_k^h) \leq \left(\frac{h \log (9m)}{\log \phi}+3\right)(4h+1).$$
Let $k_0$ be sufficiently large such that
\begin{equation}\label{k0up}
  \left(\frac{h \log (9m)}{\log \phi}+3\right)(4h+1)<\varepsilon \left(k_0-\frac{2\log m}{\log \phi}+O(1)\right)
\end{equation}
and set $m=\phi^\gamma$. Then for any $\gamma$ sufficiently large we find $k_0$ such that
$n_{k_0}< m\phi^{8k_0+1}$ satisfies $$  \frac{s_F(n_{k_0}^h)}{s_F(n_{k_0})}<\varepsilon.$$
By a direct calculation one can check that each $k=k_0$ with~(\ref{k0up}) also satisfies~(\ref{kcondition}) provided
\begin{equation}\label{eps}
  \varepsilon<\frac{h(4h+1)}{h-2},
\end{equation}
where~(\ref{eps}) is empty for $h=2$. By construction, each distinct $m$ will give rise to a distinct $n$. We therefore have
for $\gamma$ sufficiently large,
\begin{align*}
  \alpha \;&> \frac{\gamma}{8k_0+\gamma+O(1)}> \frac{\gamma}{\frac{8}{\varepsilon}(4h+1)(h\; \frac{\log 9}{\log \phi}+h\gamma+3) +16\gamma+\gamma+O(1)}\\
  &> \frac{1}{36 h^2/\varepsilon+18}.
\end{align*}
Now, suppose $\varepsilon\geq h(4h+1)/(h-2)$.
Then we conclude
$$ \alpha > \frac{\gamma}{8k_0+\gamma+O(1)}> \frac{\gamma}{\frac{8h \log 9}{\log \phi}+8h\gamma +O(1)}>\frac{1}{8h+1}.$$
This completes the proof of Theorem~\ref{theom2}.\hfill \qed

\medskip

\textit{Proof of Theorem~\ref{theom3}:} Let $n_k=mL_{2k-1}$. With the help of~(\ref{hexp}) we see that $n_k^h$ can be written as the difference of $m^h L_{h(2k-1)}$
and a positive number that is bounded by $m^h 2^{h-1} \phi^{(h-2)(2k-1)}$. In order to have terms noninterfering we suppose
that $k$ is such that
$$m^h 2^{h-1} \phi^{(h-2)(2k-1)}<F_{h(2k-1)-1-\frac{h \log m}{\log \phi}+O(1)}<\frac{\phi^{h(2k-1)+O(1)}}{\sqrt{5}m^h},$$
or equivalently,
\begin{equation}\label{kcondition2}
  k> \frac{h\log(2m^2)}{4\log \phi}+O(1).
\end{equation}
Lemma~\ref{lucasmulti} shows that for all such $k$ we have
$$s_F(n_k)\leq \frac{h\log m}{\log \phi}+O(1).$$
On the other hand, a similar calculation as in Section~\ref{sectheom1} gives
\begin{align*}
  s_F(n_k^h) &\geq 1+\left\lfloor \frac{h(2k-1)-1-h\log m/\log \phi +O(1)}{2}\right\rfloor\\
  &\quad\qquad -\frac{\log \left(m^h 2^{h-1} \phi^{(h-2)(2k-1)}\right)}{2\log \phi}-\frac{\delta'}{2}\\
  &\geq 2k-\frac{h}{\log \phi}\left(\log m+\frac{\log 2}{2}\right)+O(1).
\end{align*}
We now choose $k_0$ in a way that
$$2 k_0-\frac{h}{\log \phi}\left(\log m+\frac{\log 2}{2}\right)+O(1)>\delta \left(\frac{h \log m}{\log \phi}+O(1)\right).$$
Observe that for any $\delta>0$ each such $k=k_0$ automatically satisfies~(\ref{kcondition2}). Put $m=\phi^\gamma$. Similarly as above we get for $\gamma$ sufficiently large,
\begin{align*}
  \beta &> \frac{\gamma}{2k_0+\gamma+O(1)} > \frac{\gamma}{\delta(h\gamma+O(1))+h\gamma+\frac{h \log 2}{2 \log \phi} +\gamma+O(1)}\\
  &> \frac{1}{\delta h+h+2}.
\end{align*}
This completes the proof of Theorem~\ref{theom3}.\hfill \qed

\end{document}